\UseRawInputEncoding

\documentclass[a4paper,12pt]{article}
\usepackage{amsmath}
\usepackage[latin1]{inputenc}
\usepackage{amsmath,amsthm,amssymb}
\usepackage{enumerate}
\usepackage[pagewise]{lineno}
\usepackage[colorlinks=true]{hyperref}
\usepackage{cases}

\usepackage{cite}
\usepackage{verbatim}

\hypersetup{urlcolor=blue, citecolor=red}
\newtheorem{theorem}{Theorem}[section]

\newtheorem{lemma}{Lemma}[section]
\newtheorem{definition}{Definition}[section]

\renewcommand{\theequation}{\thesection.\arabic{equation}}

\newcommand{\Rmnum}[1]{\expandafter\@slowromancap\romannumeral #1@}
\newcommand{\eps}{\varepsilon}

\newcommand{\tme}{T_{max,\eps}}
\newcommand{\ueps}{u_\eps}
\newcommand{\veps}{v_\eps}
\newcommand{\bom}{\overline{\Omega}}
\newcommand{\Om}{\Omega}
\newcommand{\be}{\begin{equation} \label}
\newcommand{\ee}{\end{equation}}
\newcommand{\bea}{\begin{eqnarray}\label}
\newcommand{\eea}{\end{eqnarray}}
\newcommand{\bas}{\begin{eqnarray*}}
\newcommand{\eas}{\end{eqnarray*}}
\newcommand{\bit}{\begin{itemize}}
\newcommand{\eit}{\end{itemize}}
\newcommand{\nn}{\nonumber}
\newcommand{\hs}{\hspace*}
\newcommand{\io}{\int_\Omega}
\newcommand{\na}{\nabla}
\newcommand{\Del}{\Delta}
\newcommand{\wto}{\rightharpoonup}
\newcommand{\pa}{\partial}
\newcommand{\hra}{\hookrightarrow}
\newcommand{\abs}{\\[5pt]}
\newcommand{\vp}{\varphi}
\makeatother
\title{Boundedness in a taxis-consumption system involving signal-dependent motilities
 and concurrent enhancement of density-determined diffusion and cross-diffusion
}
\author{Genglin Li\footnote{E-mail address: 1182028@mail.dhu.edu.cn}\\
{\small College of Information Science and Technology, }\\
{\small Donghua University, Shanghai 201620, P.R.~China}
\and
Liangchen Wang\footnote{E-mail address: wanglc@cqupt.edu.cn}\\
{\small School of Science, Chongqing University of Posts and Telecommunications,}\\
{\small Chongqing 400065, P.R.~China}}

\date{}

\begin{document}
\baselineskip20pt \maketitle
\maketitle
\renewcommand{\theequation}{\arabic{section}.\arabic{equation}}
\catcode`@=11 \@addtoreset{equation}{section} \catcode`@=12

%The abstract of your paper
\begin{abstract}
\noindent
\begin{abstract}
\noindent
This paper is concerned with the migration-consumption taxis system involving signal-dependent motilities
\bas
 	\left\{ \begin{array}{l}	
 	u_t = \Del \big(u^m\phi(v)\big), \\[1mm]
 	v_t = \Delta v -uv,
 	\end{array} \right.
 	\qquad \qquad (\star)
\eas
in smoothly bounded domains $\Om\subset \mathbb{R}^n$, where $m>1$ and $n\ge2$.
It is shown that if $\phi\in C^3([0,\infty))$ is strictly positive on $[0,\infty)$,
for all suitably regular initial data an associated no-flux type initial-boundary value problem
possesses a globally defined bounded weak solution, provided $m>\frac{n}{2}$, which is consistent with the restriction imposed on $m$ in
corresponding signal production counterparts of $(\star)$ so as to establish the similar result.\abs
\end{abstract}
\end{abstract}
\section{Introduction}
Among an abundance of chemotaxis models that have generated active interest in recent years is the parabolic system
\be{01}
 	\left\{ \begin{array}{l}	
 	u_t = \Del \big(u^m\phi(v)\big), \\[1mm]
 	v_t = \Delta v + f(u,v),
 	\end{array} \right.
\ee
proposed to model the collective behavior of bacterial populations underlain by local sensing mechanisms,
and viewed as a particular case of the classical Keller-Segel system (\cite{fu}, \cite{liu}, \cite{KS1}, \cite{KS3}),
in comparison to which  the migration
operator in (\ref{01}) features a special link between the cell diffusivity and the chemotactic
sensitivity. \abs
From the perspective of mathematical analysis, considerable efforts have been made to understand how the marked link mentioned above
plays out in the context  where the motility function $\phi$ is exclusively influenced by chemical signals that are produced by cells
themselves, thus with $f$ taken to be $u-v$ in (\ref{01}). In particular, in the case when the random diffusive movement of cells is partially driven by the
Brownian-type diffusion regarding cells, namely, $m=1$, for a myriad of motility functions $\phi$, the literature has afforded
profound insights into the global solvability of (\ref{01}) (\cite{fujie_senba}, \cite{ahn_yoon}, \cite{fujie_jiang_ACAP2021},\cite{jiang_laurencot},
\cite{taowin_M3AS}, \cite{burger}, \cite{desvillettes}), and some close  variants
 (\cite{fujie_jiang_JDE2020},\cite{LWL},\cite{wenbin_lv_EECT},\cite{wenbin_lv_PROCA}, \cite{jiang_arxiv},\cite{liu_xu},\cite{jin_kim_wang},\cite{wang_wang});
on the other hand, though, this type of diffusion mechanisms seems limited in terms of blow-up suppression, as suggested by numerical findings (\cite{desvillettes}), compellingly evidenced by the identification of the possibility of infinite-time blow-up solutions
(\cite{fujie_jiang_CVPDE}, \cite{fujie_senba}, \cite{jin_wang}, \cite{ahn_yoon}, \cite{DLTW}),
as well as by lack of both solution theory
of (\ref{01}) with certain types of $\phi$ involving signal-dependent degeneracies, and boundedness results for  large-data
solutions even within the multi-dimensional frameworks with non-degenerate cell diffusion (\cite{taowin_M3AS}). In pursuit of delving further into
the potential blow-up-inhibiting effects of (\ref{01}) with signal production mechanisms, the enhanced cell diffusion process
resulting from an increase of $m$ was considered in \cite{win_nonlinearity}, particularly the extent to which the stabilizing effects exerted by such
a process may overbalance the simultaneously increased cross-diffusion influences, thereby leading to the global existence
and even the boundedness of solutions to (\ref{01}). In \cite{yifu_wang}, asymptotic behavior of solutions to a close relative
of (\ref{01}) with $m>1$ was investigated. \abs
The development of research related to the analysis of (\ref{01}) has gone beyond that to include a somewhat different scenario
where the emphasis is placed on the mechanism of what is called signal consumption as opposed to that of said signal production,
resulting in the chemotaxis-consumption versions of (\ref{01}), as given by
\be{00}
	\left\{ \begin{array}{l}	
	u_t = \Del \big(u^m\phi(v)\big), \\[1mm]
	v_t = \Delta v -uv.
	\end{array} \right.
\ee
A certain amount of research has centered on (\ref{00}) with $m=1$, revealing the chemotaxis-consumption interplay of the  type
in (\ref{00}) having a strong tendency to enforce  global solvability---and even bounded solutions in some situations
(\cite{li_zhao_ZAMP} ,\cite{liwin2},\cite{liwin3}).
In fact, it is shown that, under the assumptions that $\phi$ is sufficiently smooth and strictly positive on $[0,\infty)$,
even for initial data of very little regularity, (\ref{00}) possesses a global very weak solution (\cite{liwin2}),
this followed by the result
obtained in \cite{liwin3}, asserting that when imposed further regularity on $\phi$ and the initial data, in this case $\phi\in C^3([0,\infty))$ and
$(u|_{t=0},v|_{t=0})\in (W^{1,\infty}(\Om))^2$, the global very weak solution proves classical and bounded when $n\le2$,
whereas in three and higher-dimensional settings, global weak solutions of (\ref{00}) are found to exist, and especially when $n=3$,
each of these solutions becomes eventually smooth and bounded,
 regardless of the size of initial data, which however is necessary to achieve the analogous result for the chemotaxis-production
counterparts of (\ref{00}) (\cite{win_nonlinearity});
 that tendency has been further bolstered by some recent studies that keep focus on the cases with degenerate motilities in (\ref{00}),
among other things, confirming the global existence of weak and generalized solutions for a large class of $\phi$ with various decay behavior
 (\cite{win_sig_dep_mot_cons_general},\cite{W1}, \cite{win_sig_dep_mot_cons_2}, \cite{win_sig_dep_mot_cons_largetime}). \abs
{\bf Main results.} \quad
As discussed above, the picture of solution properties of (\ref{00}) with $m=1$ seems to have become clearer, and yet the information on the case
$m>1$ appears far from complete. It is, therefore, the main objective of this paper to address the query as to how far the adoption
of cell diffusion enhancement, accordingly followed by the simultaneous enhancement of the diffusion and the cross-diffusion
both implicitly included in the first equation in (\ref{00}), will lead to the establishment of bounded global solution theory for the taxis-consumption
system (\ref{00}) with strictly positive $\phi$. Bearing this in mind, we consider the initial-boundary value problem given by
\begin{eqnarray}\label{Equ(1.1)}
\left\{
\begin{array}{llll}
u_t=\Delta (u^m\phi(v)),\quad &x\in \Omega,\quad t>0,\\
v_t=\Delta v-uv,\quad &x\in\Omega,\quad t>0\\
\frac{\partial (u^m\phi(v))}{\partial \nu}=\frac{\partial v}{\partial \nu}=0,\quad &x\in\partial\Omega,\,\, t>0,\\
u(x,0)=u_0(x),\,\,\,v(x,0)=v_0(x),\quad &x\in\Omega,\\
\end{array}
\right.
\end{eqnarray}
in a smoothly bounded domain $\Omega \subset \mathbb{R}^n(n\geq2)$. \abs
In this setting, it turns out that
the restriction on $m$ identified in this paper, which is sufficient for the construction of global bounded solutions to (\ref{Equ(1.1)}),
is parallel to that made in \cite{win_nonlinearity}, and our main result thus reads as follows.
\begin{theorem}\label{result1.1}
Let $\Omega\subset \mathbb{R}^n$ ($n\geq 2$) be a bounded domain with smooth boundary, and suppose that
 \be{phi}
	\phi\in C^3([0,\infty))
	\mbox{\quad is such that \quad } \phi(\xi)>0 \mbox{ for all } \xi\ge 0,
  \ee
and that
$$m>\frac{n}{2}.$$
Then for all the initial data $(u_0,v_0)$  fulfilling that

\begin{equation}\label{Equ(1.2)}
\left\{
\begin{array}{llll}
	u_0\in W^{1,\infty}(\Om) \mbox{ is non-negative in } \bom, \\
     v_0\in W^{1,\infty}(\Omega) \mbox{ is non-negative in } \bom,
\end{array}
\right.
\end{equation}
the problem (\ref{Equ(1.1)}) admits at least one global weak solution $(u,v)$
in the sense of Definition \ref{result0.1} below, which satisfies that there exists $C>0$ in such a way that
\begin{equation*}
\|u(\cdot,t)\|_{L^{\infty}(\Omega)}+\|v(\cdot,t)\|_{W^{1,\infty}(\Omega)}\leq C \quad \text{for all } t>0.
\end{equation*}
\end{theorem}

{\bf Main ideas.} \quad
In light of the potential  degeneracies included in the first equation of (\ref{Equ(1.1)}), our analysis needs to be
based on establishing estimates for the solutions $(\ueps,\veps)$ to suitably regularized variants (\ref{Equ(2.1)})
of (\ref{Equ(1.1)}). We take the first step to derive a spatio-temporal a priori bound  for $\ueps^{m+1}$
through a duality-based argument (Lemma \ref{result3.1}), and as a result of this, an $L^q$ bound of $\veps$
for some $q>n$ will be achieved under the pivotal assumption $m>\frac{n}{2}$ (Lemma \ref{result3.2}),
paving the way for the derivation of $\|\ueps\|_{L^p(\Om)}$ for any $p\in(1,\infty)$ (Lemma \ref{result3.4}),
whence the uniform boundedness of $\ueps$ and $\na \veps$ could be conveniently derived (Lemma \ref{result3.5}).
The established uniform bounds of $\|\ueps\|_{L^{\infty}(\Om)}$ and $\|\veps\|_{W^{1,\infty}(\Om)}$
 have several applications, among them in implying some higher order regularity features
of $\ueps$ (Lemmata \ref{result3.7} and \ref{result3.8}) and in establishing H\"{o}lder estimates for $\veps$ (Lemma \ref{result3.9}).
In view of these estimates, a limit function $(u,v)$ can be identified through a limit procedure and is shown to be a global
weak solution to (\ref{Equ(1.1)}) (Lemma \ref{result3.10}).
\section{Preliminaries}
In view of the fact that cell diffusion in (\ref{Equ(1.1)}) possibly including degeneracies as $m$ increases indicates that classical solutions
might not be expected, we start with specifying the solution concept pursued in this paper.
\begin{definition}\label{result0.1}
Let $m>1$ and $\phi\in W^{1,\infty}(\Om)$ be non-negative, and suppose that $0\le u_0\in L^1(\Om)$ and $0\le v_0\in L^1(\Om)$.
A pair ($u, v$) of non-negative functions
\begin{equation}\label{Equ(1.5)}
\left\{\begin{array}{l}
u \in L_{l o c}^1(\bar{\Omega} \times[0, \infty)) \quad \text { and } \\
v \in  L_{l o c}^{\infty}(\bar{\Omega}\times[0,\infty)) \cap L_{l o c}^{1}\left([0, \infty) ; W^{1,1}(\Omega)\right)
\end{array}\right.
\end{equation}
will be called a global weak solution of (\ref{Equ(1.1)}), if
 \be{w1}
  u^m\in L_{l o c}^1([0, \infty); W^{1,1}(\Om))
 \qquad\mbox{and}\qquad
  u^m\na v \in L_{l o c}^1 (\bom \times [0,\infty)),
 \ee
 and if
\begin{equation}\label{Equ(1.6)}
-\int_{0}^{\infty} \int_{\Omega} u \varphi_{t}-\int_{\Omega} u_0 \varphi (\cdot,0)
= -\int_{0}^{\infty} \int_{\Omega} \phi(v) \na u^m\cdot \na \vp
   -\int_{0}^{\infty} \io \phi'(v) u^m \na v \cdot \na \varphi
\end{equation}
as well as
\begin{equation}\label{Equ(1.7)}
-\int_{0}^{\infty} \int_{\Omega} v \varphi_{t}-\int_{\Omega} v_{0} \varphi(\cdot, 0)=-\int_{0}^{\infty} \int_{\Omega} \nabla v \cdot \nabla \varphi-\int_{0}^{\infty} \int_{\Omega} u v \varphi
\end{equation}
hold for all $\varphi \in C_{0}^{\infty}(\bar{\Omega} \times[0, \infty))$.
\end{definition}
In order to gear toward the construction of the solutions defined above through approximation, we consider the regularized
variants of (\ref{Equ(1.1)})
\begin{eqnarray}\label{Equ(2.1)}
\left\{
\begin{array}{llll}
u_{\varepsilon t}=\Delta ((u_{\varepsilon}+\varepsilon)^m\phi(v_\varepsilon)),\quad &x\in \Omega,\quad t>0,\\
v_{\varepsilon t}=\Delta v_\varepsilon-\frac{\ueps\veps}{1+\varepsilon\ueps},\quad &x\in\Omega,\quad t>0,\\
\frac{\partial u_\varepsilon}{\partial \nu}=\frac{\partial v_\varepsilon}{\partial \nu}=0,\quad &x\in\partial\Omega,\quad t>0,\\
\ueps(x,0)=u_0(x),\,\,\,\veps(x,0)=v_0(x),\quad &x\in\Omega,\\
\end{array}
\right.
\end{eqnarray}
for $\eps\in(0,1)$, which are indeed globally solvable in the classical sense by means of the well-established parabolic theory
in \cite{A}.
\begin{lemma}\label{result2.1}
Let $\Omega\subset \mathbb{R}^{n}(n\geq2)$ be a bounded domain with smooth boundary and $m>1$, and suppose that
(\ref{phi}) and (\ref{Equ(1.2)}) hold.  Then for each $\eps\in(0,1)$, there exist
\be{lc}
\left\{
\begin{array}{llll}
u_\varepsilon\in C^0(\overline{\Omega}\times[0,\infty))\cap C^{2,1}(\overline{\Omega}\times(0,\infty)),\\
v_{\varepsilon} \in \bigcap_{q>n} C^{0}\left([0, \infty) ; W^{1, q}(\Omega)\right) \cap C^{2,1}(\bar{\Omega} \times(0, \infty))
\end{array}
\right.
\ee
such that $\ueps,\veps\ge0$ in $\bom\times[0,\infty)$,  as well as that $(\ueps,\veps)$ solves (\ref{Equ(1.1)}) classically  in $\Omega\times(0,\infty)$.
Moreover, the following properties hold
\begin{equation}\label{Equ(2.2)}
\|u_\varepsilon(\cdot,t)\|_{L^1(\Omega)}=\|u_{0}\|_{L^1(\Omega)}
	\qquad \mbox{for all $t>0$ and } \eps\in (0,1)
\end{equation}
and
\begin{equation}\label{Equ(2.3)}
\|v_\varepsilon(\cdot,t)\|_{L^{\infty}}
\leq\|v_{0}\|_{L^{\infty}(\Omega)}
	\qquad \mbox{for all $t>0$ and } \eps\in (0,1).
\end{equation}
\end{lemma}

\noindent{\bf{Proof.}} By means of  (\ref{phi}) and (\ref{Equ(1.2)}), for each $\eps\in(0,1)$,
the standard parabolic theory (\cite{amann_2000}), applicable to (\ref{Equ(2.1)}), provides $\tme\in(0,\infty]$ and at least
one classical solution satisfying (\ref{lc}) which is such that
\be{tme}
  	\mbox{if $\tme<\infty$}, \quad then \quad
  	\limsup_{t\nearrow \tme} \|\ueps(\cdot,t)\|_{L^{\infty}(\Omega)} =\infty.
  \ee
Moreover, it follows from (\ref{Equ(1.2)}) and a comparison argument that $\ueps\ge0$ and $\veps\ge0$ in $\bom\times[0,\tme)$.
Next, integrating the first equation of (\ref{Equ(2.1)}), we readily obtain
\begin{equation}\label{Equ(2.2)}
\|u_\varepsilon(\cdot,t)\|_{L^1(\Omega)}=\|u_{0}\|_{L^1(\Omega)}
	\qquad \mbox{for all $t\in(0,\tme)$ and } \eps\in (0,1),
\end{equation}
and the maximum principle, applied to the second equation of (\ref{Equ(2.1)}), implies that
\begin{equation}\label{Equ(2.3)}
\|v_\varepsilon(\cdot,t)\|_{L^{\infty}}
\leq\|v_{0}\|_{L^{\infty}(\Omega)}
	\qquad \mbox{for all $t\in(0,\tme)$ and } \eps\in (0,1).
\end{equation}
It remains to prove $\tme=\infty$ for all $\eps\in(0,1)$, which is approached through a contradiction argument
that is analogous to that given in the proof of \cite[Lemma 2.2]{liwin2}. The argument proceeds as follows:
first assuming $\tme<\infty$ for some $\eps\in(0,1)$, we use the relation
$\Big\|\frac{\ueps(\cdot,t)\veps(\cdot,t)}{1+\eps\ueps(\cdot,t)}\Big\|_{L^{\infty}(\Om)}
  \le\frac{\|v_{0}\|_{L^{\infty}(\Om)}}{\eps}$ for all $t\in (0,\tme)$
  and standard parabolic regularity theory to find that there exists $c_1(\eps)>0$ such that
   \bas
    	\|\na \veps(\cdot,t)\|_{L^{\infty}(\Om)} \le c_1(\eps)
    	\qquad \mbox{for all $t\in (0,\tme)$},
    \eas
and then with the help of this, the relations $\ueps\le\ueps+\eps\le\ueps+1$ as well as (\ref{phi}) and (\ref{Equ(2.3)}),
a standard $L^p$ testing procedure established for any $p>1$ will end up leading to the existence of some $c_2(p, \eps)>0$
fulfilling that
  \bas
   \|\ueps(\cdot,t)\|_{L^p(\Om)} \le c_2(p,\eps)
   	\qquad \mbox{for all $t\in (0,\tme)$},
  \eas
which in conjunction with a Moser-type iteration (\cite[Lemma A.1]{taowin_subcrit}) enables us to find $c_3(\eps)>0$ such that
\bas
 \|\ueps(\cdot,t)\|_{L^{\infty}(\Om)}\le c_3(\eps)
 \qquad \mbox{for all $t\in (0,\tme)$},
\eas
this clearly contradicting (\ref{tme}) and finally entailing $\tme=\infty$ for all $\eps\in(0,1)$.
\qed

In what follows, we shall fix $u_0$ and $v_0$ satisfying (\ref{Equ(1.2)}), and suppose that $\phi$ is valid with (\ref{phi}).
We begin to investigate some further basic properties of approximate solutions $(\ueps,\veps)$.

\begin{lemma}\label{result2.2}
Let $m>1$. Then
\begin{equation}\label{Equ(3.21)}
\begin{aligned}
&\int_{0}^{\infty} \int_{\Omega}\left|\nabla v_{\varepsilon}\right|^{2} \leq \frac{1}{2} \int_{\Omega} v_{0}^{2}
\qquad \mbox{for all } \eps\in (0,1),
\end{aligned}
\end{equation}
and, moreover, there exist $\kappa_1\in(0,\infty), \kappa_2\in(0,\infty)$ and $\kappa_3\in[0,\infty)$ such that for all $\eps\in(0,1)$
\begin{equation}\label{Equ(2.4)}
\kappa_1\leq \phi(v_\varepsilon(x,t))\leq \kappa_2
	\qquad \mbox{for all $x\in\Om$ and $t>0$}
\end{equation}
and
\begin{equation}\label{Equ(2.5)}
|\phi'(v_\varepsilon(x,t))|\leq \kappa_3
		\qquad \mbox{for all $x\in\Om$ and $t>0$}.
\end{equation}
\end{lemma}
\noindent{\bf{Proof.}}
Multiplying the second equation of (\ref{Equ(2.1)}) by $\veps$ and integrating over $t>0$ we have
\bas
\frac{1}{2} \int_{\Omega} v_{\varepsilon}^{2}(\cdot,t)
+ \int_{0}^{t}\int_{\Omega}\left|\nabla v_{\varepsilon}\right|^{2}
+ \int_{0}^{t} \int_{\Omega} \frac{u_{\varepsilon} v_{\varepsilon}^2}{1+\eps\ueps}
\le \frac{1}{2}\io v_0^2
\eas
for all $\eps\in(0,1)$, and therefore, we obtain (\ref{Equ(3.21)}).
Consequently, employing (\ref{phi}) and (\ref{Equ(2.3)}) yields
(\ref{Equ(2.4)}) and (\ref{Equ(2.5)}) with
$\kappa_1:=\mathop {\min }\limits_{0\leq s \leq\|v_0\|_{L^\infty(\Omega)} }\phi(s)>0$,
$\kappa_2:=\mathop {\max }\limits_{0\leq s \leq\|v_0\|_{L^\infty(\Omega)} }\phi(s)>0$
and
 $\kappa_3:= \mathop {\max }\limits_{0\leq s \leq\|v_0\|_{L^\infty(\Omega)} }|\phi'(s)|\ge0$. $\hfill{} \Box$

\section{Proof of Theorem \ref{result1.1} }
\subsection{ A priori estimates }

The purpose of this section is to prove the uniform boundedness of approximate solutions to (\ref{Equ(2.1)}).
In the first step, by means of a duality-based argument derived from that used in \cite{taowin_M3AS}, we can obtain a spatio-temporal $L^1$
bound for $u^{m+1}_\varepsilon$, which is of essential importance to the construction of bounded solutions.

\begin{lemma}\label{result3.1}
Assume $m>1$, then there exists $C>0$ such that
\begin{equation}\label{Equ(3.1)}
\int_{t}^{t+1}\int_\Omega u_\varepsilon^{m+1}\leq C
	\qquad \mbox{for all $t>0$ and } \eps\in (0,1).
\end{equation}
\end{lemma}
\noindent{\bf{Proof.}} Let $A$ be the self-adjoint realization of $-\Delta+1$ under homogeneous Neumann boundary conditions in $L^{2}(\Omega)$ (see more details in \cite{S}). Hence, we can rewrite the first equation of (\ref{Equ(2.1)}) as follows
 \bas
 u_{\varepsilon t}+A[(u_\varepsilon+\varepsilon)^{m}\phi(v_\varepsilon)]=(u_\varepsilon+\varepsilon)^{m}\phi(v_\varepsilon)
 \quad \mbox{for all $t>0$ and }\eps\in(0,1).
 \eas
Testing this equation against $A^{-1}(u_\varepsilon+\varepsilon)$, and using $A$ being self-adjoint as well as Young's inequality, we can find $c_1>0$ such that
for all $t>0$ and $\eps\in(0,1)$
\bas
\dfrac{1}{2}\dfrac{d}{dt}{\int _{\Omega}|A^{-\frac{1}{2}}(u_\varepsilon+\varepsilon)|^2}
+\int_\Omega (u_\varepsilon+\varepsilon)^{m+1}\phi(v_\varepsilon)
&=&\int_\Omega (u_\varepsilon+\varepsilon)^{m}\phi(v_\varepsilon)A^{-1}(u_\varepsilon+\varepsilon)\\
&\le&\frac{1}{2}\int_\Omega (u_\varepsilon+\varepsilon)^{m+1}\phi(v_\varepsilon)\\
  & &+c_1\int_\Omega \phi(v_\varepsilon)|A^{-1}(u_\varepsilon+\varepsilon)|^{m+1},
\eas
hence by (\ref{Equ(2.4)}) and (\ref{Equ(2.5)}), we have
\begin{equation}\label{Equ(3.2)}
\begin{aligned}
\dfrac{d}{dt}{\int _{\Omega}|A^{-\frac{1}{2}}(u_\varepsilon+\varepsilon)|^2}+\kappa_1\int_\Omega (u_\varepsilon+\varepsilon)^{m+1}&+\int _{\Omega}|A^{-\frac{1}{2}}(u_\varepsilon+\varepsilon)|^2\\
\leq&2c_1\kappa_2\int_\Omega|A^{-1}(u_\varepsilon+\varepsilon)|^{m+1}+\int_{\Omega}|A^{-\frac{1}{2}}(u_\varepsilon+\varepsilon)|^2
\end{aligned}
\end{equation}
for all $t>0$ and $\eps\in(0,1)$. Next we will show that the integrals on the right-hand side can be controlled by $\int_\Omega (u_\varepsilon+\varepsilon)^{m+1}$.
To this end, letting $q>1$ be such that $\frac{n(m+1)}{n+2m+2}\le q< m+1$, we derive from the left inequality  that the embedding
$ W^{2,q}(\Omega)\hookrightarrow L^{m+1}(\Omega)$ is continuous, which combined with standard elliptic regularity in $L^{m+1}(\Omega)$ (\cite{GT}) warrants
the existence of  $c_2>0$ and $c_3 > 0$ such that
\bas
\|\psi\|_{L^{m+1}(\Omega)} \leq c_2\|\psi\|_{W^{2,q}(\Omega)}\leq c_3\|A\psi\|_{L^{q}(\Omega)}
\qquad\mbox{for all $\psi \in W^{2,q}(\Omega)$ fulfilling $\frac{\partial \psi}{\partial \nu}=0$ on $\partial\Omega$}.
\eas
Taking into account this inequality along with
 $1<q< m+1$, we use Young's inequality and (\ref{Equ(2.2)}) to find $c_4>0$ and $c_5>0$ in such a way that
\begin{equation}\label{Equ(3.4)}
\begin{aligned}
2c_1\kappa_2\|A^{-1}(u_\varepsilon+\varepsilon)\|^{m+1}_{L^{m+1}(\Omega)}
\leq& 2c_1\kappa_2c^{m+1}_3\|u_\varepsilon+\varepsilon\|^{m+1}_{L^{q}(\Omega)}\\
\leq& 2c_1\kappa_2c^{m+1}_3\|u_\varepsilon+\varepsilon\|^{(m+1)\theta}_{L^{m+1}(\Omega)}\|u_\varepsilon+\varepsilon\|^{(m+1)(1-\theta)}_{L^1(\Omega)}\\
\leq&\frac{\kappa_1}{4}\|u_\varepsilon+\varepsilon\|^{m+1}_{L^{m+1}(\Omega)}+c_4\|u_\varepsilon+\varepsilon\|^{m+1}_{L^{1}(\Omega)}\\
\leq& \frac{\kappa_1}{4}\|u_\varepsilon+\varepsilon\|^{m+1}_{L^{m+1}(\Omega)}+c_5
\quad \mbox{for all $t>0$ and }\eps\in(0,1)
\end{aligned}
\end{equation}
with $\theta:=\frac{1-\frac{1}{q}}{1-\frac{1}{m+1}}\in(0,1)$.
Now picking $p>1$ such that $p\ge\frac{2n}{n+2}$, and setting $p'=\frac{p}{p-1}$, we thus infer that $(n-2p)p'\le np$,
which implies the continuity of the embedding
$ W^{2,p}(\Omega)\hookrightarrow L^{p'}(\Omega)$. In view of this,  standard elliptic regularity provides positive constants
$c_6$ and $c_7$ satisfying that
\begin{equation}\label{Equ(3.6)}
\begin{aligned}
\|\psi\|_{L^{p'}(\Omega)} \leq c_6\|\psi\|_{W^{2,p}(\Omega)}\leq c_7\|A\psi\|_{L^{p}(\Omega)}
\end{aligned}
\end{equation}
for all $\psi \in W^{2,p}(\Omega)$ such that $\frac{\partial \psi}{\partial \nu}=0$ on $\partial\Omega$. Since $m>1$ implies $m+1>2>\frac{2n}{n+2}$,
by Young's inequality, we may find $c_8>0$ fulfilling that
\begin{equation}\label{Equ(3.7)}
\begin{aligned}
\int_{\Omega}|A^{-\frac{1}{2}}(u_\varepsilon+\varepsilon)|^2=\int_{\Omega}(u_\varepsilon+\varepsilon)A^{-1}(u_\varepsilon+\varepsilon)
&\leq \|u_\varepsilon+\varepsilon\|_{L^{m+1}(\Omega)}\|A^{-1}(u_\varepsilon+\varepsilon)\|_{L^\frac{m+1}{m}(\Omega)}\\
&\leq c_7\|u_\varepsilon+\varepsilon\|^2_{L^{m+1}(\Omega)}\\
& \leq \frac{\kappa_1}{4}\|u_\varepsilon+\varepsilon\|^{m+1}_{L^{m+1}(\Omega)} +c_8
\end{aligned}
\end{equation}
for all $t>0$ and $\eps\in(0,1)$. Plugging (\ref{Equ(3.4)}) and (\ref{Equ(3.7)}) back into (\ref{Equ(3.2)}) enables us to find $c_9>0$ such that
\begin{equation}\label{Equ(3.9)}
\begin{aligned}
\dfrac{d}{dt}{\int _{\Omega}|A^{-\frac{1}{2}}(u_\varepsilon+\varepsilon)|^2}+\frac{\kappa_1}{2}\int_\Omega (u_\varepsilon+\varepsilon)^{m+1}+\int _{\Omega}|A^{-\frac{1}{2}}(u_\varepsilon+\varepsilon)|^2\leq c_9
\end{aligned}
\end{equation}
for all $t>0$ and $\eps\in(0,1)$. By means of an ODE comparison argument, we have
\be{311}
\int_{\Omega}\left|A^{-\frac{1}{2}}\left(u_{\varepsilon}+\eps\right)(\cdot,t)\right|^{2}
\leqslant c_{10}:=\max \left\{ \sup_{\varepsilon\in (0,1)}\int_{\Omega}\left|A^{-\frac{1}{2}}\left(u_{0}+\eps\right)\right|^{2}, \frac{2c_{9}}{\kappa_{1}}\right\}
\ee
for all $t>0$ and $\eps\in(0,1)$, with $c_{10}$ being finite as guaranteed by the continuity of $A^{-\frac{1}{2}}$
and (\ref{Equ(1.2)}), and consequently, by (\ref{311}) we conclude (\ref{result3.1}) from a direct integration in (\ref{Equ(3.9)}).
\qed

Relying on the fundamental assumption $m>\frac{n}{2}$, we can establish bounds for $\na v_{\eps}$
with respect to the norm in $L^{q}(\Om)$ for some $q>n$ on the basis of the information derived in Lemma \ref{result3.1}.
\begin{lemma}\label{result3.2}
If $m>\frac{n}{2}$, then there exist $q>n$ and $C>0$ such that
\begin{equation*}
\|\nabla v_{\eps}(\cdot,t)\|_{L^{q}(\Omega)}\leq C
\quad \qquad \mbox{for all $t>0$ and } \eps\in (0,1).
\end{equation*}
\end{lemma}
\noindent{\bf{Proof.}} Since, by assumption, $m>\frac{n}{2}$, we have
$$\Big\{-\frac{1}{2}-\frac{n}{2}\Big(\frac{1}{m+1}-\frac{1}{n}\Big)\Big\}\cdot \frac{m+1}{m}>-1.$$
Thus we may fix $q>n$ sufficiently close to $n$ so that
\begin{equation}\label{Equ(3.132)}
\Big\{-\frac{1}{2}-\frac{n}{2}\Big(\frac{1}{m+1}-\frac{1}{q}\Big)\Big\}\cdot \frac{m+1}{m}>-1.
\end{equation}
In addition, from $m>1$ it follows that
\bea{3130}
-\frac{1}{2}\cdot \frac{m+1}{m}>-1.
\eea
 Using (\ref{Equ(3.132)}) in conjunction with (\ref{3130}) allows us to pick $p\le\min\{q,m+1\}$ with the properties that $p\ge1$ and
\be{3131}
 \Big\{-\frac{1}{2}-\frac{n}{2}\Big(\frac{1}{p}-\frac{1}{q}\Big)\Big\}\cdot \frac{m+1}{m}>-1.
\ee
Then we invoke the smoothing properties of Neumann heat semigroup $\left(e^{t \Delta}\right)_{t \geq 0}$ on $\Omega$ (\cite[Lemma 1.3]{win_JDE}) to find positive constants $c_1$, $c_{2}$ and $c_{3}$ such that
\begin{align}
\|\nabla e^{t\Delta}\psi\|_{L^{q}(\Omega)}\leq c_1\|\psi\|_{W^{1,\infty}(\Omega)}
&\quad \text{for all } t>0 \,\,\text{and } \psi\in W^{1,\infty}(\Omega),\label{Equ(3.10)}\\
\|\nabla e^{\Delta}\psi\|_{L^{q}(\Omega)}\leq c_2\|\psi\|_{L^{q}(\Omega)}
&\quad \text{for all } \psi\in L^{q}(\Omega),\label{Equ(3.11)}
\end{align}
as well as
\be{Equ(3.12)}
\|\nabla e^{t\Delta}\psi\|_{L^{q}(\Omega)}
\leq c_3t^{-\frac{1}{2}-\frac{n}{2}(\frac{1}{p}-\frac{1}{q})}\|\psi\|_{L^{p}(\Omega)}
\quad \text{for all } t>0 \,\,\text{and } \psi\in C^0(\bom).
\ee
By virtue of the variation-of-constants formula associated with the second equation in (\ref{Equ(2.1)}), we see that
\begin{equation}\label{Equ(3.14)}
\begin{aligned}
\|\nabla v_\varepsilon(\cdot,t)\|_{L^{q}(\Omega)}\leq& \|\nabla e^{(t-(t-1)_+)\Delta}v_\varepsilon(\cdot,(t-1)_+)\|_{L^{q}(\Omega)}\\
&+\int_{(t-1)_+}^{t} \bigg\|\nabla e^{(t-s)\Delta}\frac{u_\varepsilon(\cdot,s)v_\varepsilon(\cdot,s)}{1+\eps u_\varepsilon(\cdot,s)}\bigg\|_{L^{q}(\Omega)}ds\\
=:&I_1+I_2
\quad \qquad \mbox{for all $t>0$ and } \eps\in (0,1).
\end{aligned}
\end{equation}
If $t\leq1$, using (\ref{Equ(3.10)}) we have
\bea{Equ(3.15)}
I_1 &=&\|\nabla e^{t\Delta}v_\varepsilon(\cdot,0)\|_{L^{q}(\Omega)}\nn\\
     &\leq& c_1\|v_0\|_{W^{1,\infty}(\Omega)}
     \quad \qquad \mbox{for all $t>0$ and } \eps\in (0,1),
\eea
while for the case  $t>1$, (\ref{Equ(3.11)}) and (\ref{Equ(2.3)}) imply  that
\bea{Equ(3.16)}
I_1 &=&\|\nabla e^{\Delta}v_\varepsilon(\cdot,t-1)\|_{L^{q}(\Omega)}\nn\\
     &\leq& c_2\|v_\varepsilon(\cdot,t-1)\|_{L^{q}(\Omega)}\nn\\
     &\leq& c_2|\Omega|^\frac{1}{q}\|v_0\|_{L^{\infty}(\Omega)}
     \quad \qquad \mbox{for all $t>0$ and } \eps\in (0,1).
\eea
Furthermore, using (\ref{Equ(3.12)}), (\ref{Equ(2.3)}), and H\"{o}lder's inequality, we see that
\begin{equation}\label{Equ(3.17)}
\begin{aligned}
I_2\leq& c_3\int_{(t-1)_+}^{t}(t-s)^{-\frac{1}{2}-\frac{n}{2}(\frac{1}{p}-\frac{1}{q})}\|u_\varepsilon(\cdot,s)v_\varepsilon(\cdot,s)\|_{L^{p}(\Omega)}ds\\
\leq& c_3\|v_0\|_{L^{\infty}(\Omega)}|\Omega|^\frac{m+1-p}{(m+1)p}\int_{(t-1)_+}^{t}(t-s)^{-\frac{1}{2}-\frac{n}{2}(\frac{1}{p}-\frac{1}{q})}\|u_\varepsilon(\cdot,s)\|_{L^{m+1}(\Omega)}ds\\
\leq& c_3\|v_0\|_{L^{\infty}(\Omega)}|\Omega|^\frac{m+1-p}{(m+1)p}\bigg\{\int_{(t-1)_+}^{t}(t-s)^{\big[-\frac{1}{2}-\frac{n}{2}\big(\frac{1}{p}-\frac{1}{q}\big)\big]
\cdot{\frac{m+1}{m}}}ds\bigg\}^{\frac{m}{m+1}}\\
&\cdot\bigg\{\int_{(t-1)_+}^{t}\|u_\varepsilon(\cdot,s)\|^{m+1}_{L^{m+1}(\Omega)}ds\bigg\}^{\frac{1}{m+1}}\\
\leq& c_3\|v_0\|_{L^{\infty}(\Omega)}|\Omega|^\frac{m+1-p}{(m+1)p}\bigg\{\int_{0}^{1}\sigma^{\big[-\frac{1}{2}-\frac{n}{2}\big(\frac{1}{p}-\frac{1}{q}\big)\big]
\cdot{\frac{m+1}{m}}}ds\bigg\}^{\frac{m}{m+1}}\\
&\cdot\bigg\{\int_{(t-1)_+}^{t}\|u_\varepsilon(\cdot,s)\|^{m+1}_{L^{m+1}(\Omega)}ds\bigg\}^{\frac{1}{m+1}}\\
\leq& c_4
\quad \qquad \mbox{for all $t>0$ and } \eps\in (0,1)
\end{aligned}
\end{equation}
with $$c_4:= c_3\|v_0\|_{L^{\infty}(\Omega)}|\Omega|^\frac{m+1-p}{(m+1)p}\Big\{\int_{0}^{1}\sigma^{[-\frac{1}{2}-\frac{n}{2}(\frac{1}{p}-\frac{1}{q})]
\cdot{\frac{m+1}{m}}}ds\Big\}^{\frac{m}{m+1}} \cdot
\Big\{\int_{(t-1)_+}^{t}\|u_\varepsilon(\cdot,s)\|^{m+1}_{L^{m+1}(\Omega)}ds\Big\}^{\frac{1}{m+1}}$$ being finite thanks to (\ref{Equ(3.1)}) and (\ref{3131}).
Collecting (\ref{Equ(3.14)})-(\ref{Equ(3.17)}), we can thereupon prove the result. \qed

The preceding lemma can now be used to derive an $L^p$ bound for the first solution component, thereby eventually leading to its uniform boundedness.
This is accomplished by the following lemma, which will be also used in the derivation of regularity feature of time derivative of $u_{\varepsilon}$.
\begin{lemma}\label{result3.3}
Let $p>0$ and $\varphi \in C^{\infty}(\Om)$. Then for all $t>0$ and $\eps\in (0,1)$
\bas
\frac{1}{p} \int_{\Omega} \frac{\partial}{\partial t}\left(u_{\varepsilon}+\varepsilon\right)^{p}\varphi
&=&-m(p-1) \int_{\Omega}\left(u_{\varepsilon}+\varepsilon\right)^{m+p-3}\phi(v_{\varepsilon})\left|\nabla u_{\varepsilon}\right|^{2} \varphi \\
& &-(p-1) \int_{\Omega}\left(u_{\varepsilon}+\varepsilon\right)^{m+p-2} \phi'(v_{\varepsilon})\varphi \nabla u_{\varepsilon} \cdot \nabla v_{\varepsilon}\\
& &-m\int_{\Omega}\left(u_{\varepsilon}+\varepsilon\right)^{m+p-2} \phi(v_{\varepsilon}) \nabla u_{\varepsilon} \cdot \nabla \varphi \\
& &-\int_{\Omega} \left(u_{\varepsilon}+\varepsilon\right)^{m+p-1} \phi'(v_{\varepsilon}) \nabla v_{\varepsilon} \cdot \nabla \varphi.
\eas
\end{lemma}
\noindent{\bf{Proof.}}  Through a direct computation we see that
\bas
& & \hs{-10mm}
\frac{1}{p} \int_{\Omega} \frac{\partial}{\partial t}\left(u_{\varepsilon}+\varepsilon\right)^{p} \cdot \varphi\\
&=&-\int_{\Omega}\bigg\{\Big\{m(u_{\varepsilon}+\varepsilon)^{m-1}\phi(v_{\varepsilon})\nabla u_{\varepsilon}
+(u_{\varepsilon}+\varepsilon)^{m} \phi'(v_{\varepsilon}) \nabla v_{\varepsilon}\Big\}\\
& & \cdot\Big\{(p-1)(u_{\varepsilon}+\varepsilon)^{p-2}\varphi\nabla u_{\varepsilon}+(u_{\varepsilon}+\varepsilon)^{p-1} \nabla \varphi\Big\}\bigg\}\\
&=&-m(p-1) \int_{\Omega}\left(u_{\varepsilon}+\varepsilon\right)^{m+p-3}\phi(v_{\varepsilon})\left|\nabla u_{\varepsilon}\right|^{2} \varphi
-\int_{\Omega} \left(u_{\varepsilon}+\varepsilon\right)^{m+p-1} \phi'(v_{\varepsilon}) \nabla v_{\varepsilon} \cdot \nabla \varphi \\
& &-m\int_{\Omega}\left(u_{\varepsilon}+\varepsilon\right)^{m+p-2} \phi(v_{\varepsilon}) \nabla u_{\varepsilon} \cdot \nabla \varphi-(p-1) \int_{\Omega}\left(u_{\varepsilon}+\varepsilon\right)^{m+p-2} \phi'(v_{\varepsilon}) \varphi\nabla u_{\varepsilon} \cdot \nabla v_{\varepsilon}
\eas
for $\varphi \in C^{\infty}(\Om)$ as well as all $t>0$ and $\eps\in(0,1)$. The lemma is thus proved. $\hfill{} \Box$

As already pointed out, the first application of Lemma \ref{result3.3} occurs when $\varphi=1$, resulting in bounds for $\ueps$ in $L^p(\Om)$ for arbitrarily
large $p$.
\begin{lemma}\label{result3.4}
Let $m>\frac{n}{2}$. Then for all  $p>1$, there exists $C(p)>0$ such that
$$
\left\|u_{\varepsilon}(\cdot, t)\right\|_{L^{p}(\Omega)} \leq C(p)
\quad \qquad \mbox{for all $t>0$ and } \eps\in (0,1).
$$
\end{lemma}
\noindent{\bf{Proof.}} Letting $\varphi=1$ in Lemma \ref{result3.3}, we then use Young's inequality to see that
\begin{equation}\label{Equ(3.180)}
\begin{aligned}
\frac{d}{dt} &\int_{\Omega}\left(u_{\varepsilon}+\varepsilon\right)^{p}+mp(p-1)\int_{\Omega}\left(u_{\varepsilon}+\varepsilon\right)^{m+p-3}\phi(v_{\varepsilon})\left|\nabla u_{\varepsilon}\right|^{2}\\
=&-p(p-1)\int_{\Omega}\left(u_{\varepsilon}+\varepsilon\right)^{m+p-2} \phi'(v_{\varepsilon}) \nabla u_{\varepsilon} \cdot \nabla v_{\varepsilon}\\
\leq&\frac{mp(p-1)}{2}\int_{\Omega}\left(u_{\varepsilon}+\varepsilon\right)^{m+p-3}\phi(v_{\varepsilon})\left|\nabla u_{\varepsilon}\right|^{2}+\frac{p(p-1)}{2m}\int_{\Omega}\left(u_{\varepsilon}+\varepsilon\right)^{m+p-1}\frac{\phi'^2(v_{\varepsilon})}{\phi(v_{\varepsilon})}\left|\nabla v_{\varepsilon}\right|^{2}
\end{aligned}
\end{equation}
for all $t>0$ and $\eps\in(0,1)$. Furthermore, recalling (\ref{Equ(2.4)}) and (\ref{Equ(2.5)}), we have
\begin{equation}\label{Equ(3.18)}
\begin{aligned}
\frac{d}{dt} &\int_{\Omega}\left(u_{\varepsilon}+\varepsilon\right)^{p}+\frac{2\kappa_1mp(p-1)}{(m+p-1)^2}\int_{\Omega}\left|\nabla\left(u_{\varepsilon}+\varepsilon\right)^\frac{m+p-1}{2} \right|^{2}+\int_{\Omega}\left(u_{\varepsilon}+\varepsilon\right)^{p}\\
\leq&\frac{\kappa^2_3p(p-1)}{2m \kappa_1}\int_{\Omega}\left(u_{\varepsilon}+\varepsilon\right)^{m+p-1}\left|\nabla v_{\varepsilon}\right|^{2}+\int_{\Omega}\left(u_{\varepsilon}+\varepsilon\right)^{p}
\end{aligned}
\end{equation}
for all $t>0$ and $\eps\in(0,1)$. Taking $q$ as provided in Lemma \ref{result3.2}, we obtain $\frac{2 q}{q-2}<\frac{2n}{(n-2)_{+}}$
and therefore the compactness of the embedding $W^{1,2}(\Omega)\hookrightarrow L^{\frac{2 q}{q-2}}(\Omega)$,
and with the help of the latter, the compactness of the embedding $W^{1,2}(\Omega)\hookrightarrow L^2(\Om)$,
and the assumption $m>\frac{n}{2}$, it follows from Lemma \ref{result3.2} and (\ref{Equ(2.2)}) that
by invoking the H\"older inequality, Young's inequality and an Ehrling-type inequality,
for all $p>1$, we can find positive constants $c_1(p)$, $c_2(p)$ and $c_3(p)$ such that
\begin{equation}\label{Equ(3.19)}
\begin{aligned}
&\frac{m\kappa^2_3p(p-1)}{2m \kappa_1}\int_{\Omega}\left(u_{\varepsilon}+\varepsilon\right)^{m+p-1}\left|\nabla v_{\varepsilon}\right|^{2}+\int_{\Omega}\left(u_{\varepsilon}+\varepsilon\right)^{p}\\
\leq&\frac{m\kappa^2_3p(p-1)}{2m \kappa_1} \left(\int_{\Omega}(u_{\varepsilon}+\varepsilon)^{\frac{(m+p-1) q}{q-2}}\right)^\frac{q-2}{q}
\left(\int_{\Omega}|\nabla v_{\varepsilon}|^{q}\right)^\frac{2}{q}+\int_{\Omega}\left(u_{\varepsilon}+\varepsilon\right)^{m+p-1}+c_1(p)\\
=& \frac{m\kappa^2_3p(p-1)}{2m \kappa_1} \left\|\left(u_{\varepsilon}+\varepsilon\right)^{\frac{m+p-1}{2}}\right\|^{2}_{L^{\frac{2 q}{q-2}}(\Omega)}\left\|\nabla v_{\varepsilon}\right\|_{L^{q}(\Omega)}^{2}+\left\|\left(u_{\varepsilon}+\varepsilon\right)^{\frac{m+p-1}{2}}\right\|^{2}_{L^{2}(\Omega)}+c_1(p)\\
\leq&c_{2}(p)\bigg\{ \left\|\left(u_{\varepsilon}+\varepsilon\right)^{\frac{m+p-1}{2}}\right\|^{2}_{L^{\frac{2 q}{q-2}}(\Omega)}+ \left\|\left(u_{\varepsilon}+\varepsilon\right)^{\frac{m+p-1}{2}}\right\|^{2}_{L^{2}(\Omega)}\bigg\}+c_1(p)\\
\leq& \frac{2\kappa_1mp(p-1)}{(m+p-1)^2}\int_{\Omega}\left|\nabla\left(u_{\varepsilon}+\varepsilon\right)^{\frac{m+p-1}{2}}\right|^{2}+c_{3}(p)
 \quad \qquad \mbox{for all $t>0$ and } \eps\in (0,1).
\end{aligned}
\end{equation}
Inserting (\ref{Equ(3.19)}) into (\ref{Equ(3.18)}), we infer that
\begin{equation*}
\begin{aligned}
\frac{d}{d t} \int_{\Omega}\left(u_{\varepsilon}+\varepsilon\right)^{p}+\int_{\Omega}\left(u_{\varepsilon}+\varepsilon\right)^{p} \leq c_{3}(p)
 \quad \qquad \mbox{for all $t>0$ and } \eps\in (0,1),
\end{aligned}
\end{equation*}
which along with a standard comparison argument implies that
$$\int_{\Omega} u_{\varepsilon}^{p}(\cdot, t) \leq \max \left\{c_{3}(p),\io (u_{0}+1)^p \right\}
 \quad \qquad \mbox{for all $t>0$ and } \eps\in (0,1).$$
Hence, we complete this proof. $\hfill{} \Box$

Thanks to the non-degenerate behavior of the form in (\ref{Equ(2.1)}), we are able to perform parabolic theory
to obtain the uniform bounds for $\ueps$ and $\na \veps$, based on Lemmata \ref{result3.3} and \ref{result3.4}.
\begin{lemma}\label{result3.5}
Let $m>\frac{n}{2}$. Then the system (\ref{Equ(2.1)}) admits a global classical solution $(u_{\varepsilon},v_{\varepsilon})$ for which $u_{\varepsilon}$ and $v_{\varepsilon}$ are global bounded in the sense that there exists a constant $K>0$ satisfying that
\begin{equation*}
\left\|u_{\varepsilon}(\cdot, t)\right\|_{L^{\infty}(\Omega)}+\left\|v_{\varepsilon}(\cdot, t)\right\|_{W^{1, \infty}(\Omega)}\leq K
\qquad \text { for all $t>0$ and $\eps\in(0,1)$}.
\end{equation*}
\end{lemma}
\noindent{\bf{Proof.}} The uniform boudedness of $\left\|v_{\varepsilon}(\cdot, t)\right\|_{W^{1, \infty}(\Omega)}$
for all $t>0$ and $\eps\in(0,1)$ is a consequence of Lemma \ref{result3.4} and the standard parabolic regularity theory
applied to the second equation of system (\ref{Equ(2.1)}) (see \cite[Lemma 4.1]{DW} for instance), and this information,
in conjunction with (\ref{Equ(2.4)}), (\ref{Equ(2.5)}) and Lemma \ref{result3.4} once more, yields the uniform boundedness of
$\left\|u_{\varepsilon}(\cdot, t)\right\|_{L^{\infty}(\Omega)}$
for all $t>0$ and $\eps\in(0,1)$ according to \cite[Lemma A.1]{taowin_subcrit}. $\hfill{} \Box$

\subsection{Further regularity properties}

As a consequence of Lemma \ref{result3.5}, some information on the gradient of $\ueps$ becomes available in the following manner.
\begin{lemma}\label{result3.7}
Let $m>\frac{n}{2}$. Then whenever $\alpha>\frac{m}{2}$, one can find $C(\alpha)>0$ such that
\begin{equation}
\int_{0}^{\infty} \int_{\Omega}\left|\nabla\left(u_{\varepsilon}+\varepsilon\right)^{\alpha}\right|^{2}  \leq C(\alpha)
\end{equation}
for all $\varepsilon \in(0,1)$.
\end{lemma}
\noindent{\bf{Proof.}} Noting that $\alpha>\frac{m}{2}$, and setting $p:= 2\alpha-m+1$, we obtain $p>1$.
 With the help of (\ref{Equ(2.4)}) and (\ref{Equ(2.5)}), it is readily seen from (\ref{Equ(3.180)}) that
\bea{Equ(3.24)}
\frac{d}{dt} \int_{\Omega}\left(u_{\varepsilon}+\varepsilon\right)^{p}
&\le& -\frac{2\kappa_1mp(p-1)}{(m+p-1)^2}\int_{\Omega}\left|\nabla\left(u_{\varepsilon}+\varepsilon\right)^\frac{m+p-1}{2} \right|^{2}\nn\\
  & & + \frac{p(p-1)\kappa^2_3}{2m\kappa_1}\int_{\Omega}\left(u_{\varepsilon}+\varepsilon\right)^{m+p-1}\left|\nabla v_{\varepsilon}\right|^{2}
\eea
for all $t>0$ and $\eps\in(0,1)$. Since Lemma \ref{result3.5} warrants the existence of $c_1>0$ such that
 $\left\|u_{\varepsilon}(\cdot, t)\right\|_{L^{\infty}(\Omega)}\le c_1$ for all $t>0$ and $\eps\in(0,1)$, we infer that
\begin{equation}\label{Equ(3.25)}
\begin{aligned}
 \frac{p(p-1)\kappa^2_3}{2m\kappa_1} \int_{\Omega}\left(u_{\varepsilon}+\varepsilon\right)^{m+p-1}\left|\nabla v_{\varepsilon}\right|^{2}
\leq \frac{p(p-1)\kappa^2_3}{2m\kappa_1}(c_1+1)^{p+m-1} \int_{\Omega}\left|\nabla v_{\varepsilon}\right|^{2}.
\end{aligned}
\end{equation}
Inserting (\ref{Equ(3.25)}) into (\ref{Equ(3.24)}) and integrating with respect to $t$ yields
\begin{equation*}
\begin{aligned}
&\int_{\Omega}\left(u_{\varepsilon}+\varepsilon\right)^{p}+\frac{2\kappa_1m p(p-1)}{(m+p-1)^{2}} \int_{0}^{t} \int_{\Omega}\left|\nabla\left(u_{\varepsilon}+\varepsilon\right)^{\frac{m+p-1}{2}}\right|^{2} \\
&\leq \frac{p(p-1)\kappa^2_3}{2m\kappa_1}(c_1+1)^{m+p-1} \int_{0}^{t} \int_{\Omega}\left|\nabla v_{\varepsilon}\right|^{2}+\int_{\Omega}\left(u_{0}+\varepsilon\right)^{p}\\
&\leq \frac{p(p-1)\kappa^2_3}{2m\kappa_1}(c_1+1)^{m+p-1}\left(\frac{1}{2} \int_{\Omega} v_{0}^{2}\right)+\int_{\Omega}\left(u_{0}+1\right)^{p}
\end{aligned}
\end{equation*}
for all $t>0$ and $\eps\in(0,1)$, where we have used (\ref{Equ(3.21)}). Recalling the definition of $p$, we thus complete the proof. $\hfill{} \Box$

In order to facilitate a compactness argument based on an Aubin-Lions lemma, we will draw on the conclusions of Lemma \ref{result3.3} and
the above estimates, thereby obtaining certain time regularity of $u_{\varepsilon}+\varepsilon$.
\begin{lemma}\label{result3.8}
Let $m>\frac{n}{2}$, and assume that $k\in \mathbb{N}$ satisfying $k>\frac{n}{2}$. Then
for any $\alpha>1$ %%$\alpha>\max \left\{\frac{m}{2}, 1\right\}$,
and for all $T \in(0, \infty)$, there exists $C(\alpha,T)>0$ such that
\begin{equation*}
\left\|\frac{\partial}{\partial t}\left(u_{\varepsilon}+\varepsilon\right)^{\alpha}\right\|_{L^{1}\left((0, T) ;\left(W_{0}^{k,2}(\Omega)\right)^{*}\right)}
\leq C(\alpha,T)
\quad \text { for all } \eps \in(0, 1).
\end{equation*}
\end{lemma}
\noindent{\bf{Proof.}} Since $k>\frac{n}{2}$ implies the continuity of the embedding $W_{0}^{k,2}(\Omega) \hookrightarrow L^{\infty}(\Omega)$,
we may find $c_1>0$ with the property that $\|\varphi\|_{L^{\infty}(\Omega)} + \|\na \varphi\|_{L^2(\Om)} \leq c_1$
for all  $\varphi \in C_0^{\infty}(\Om)$ satisfying $\|\varphi\|_{W_0^{k,2}(\Omega)}\le1$. Fixing any such $\varphi$
and using Lemma \ref{result3.3}, we find
\bea{Equ(3.26)}
\int_{\Omega} \frac{\partial}{\partial t}\left(u_{\varepsilon}+\varepsilon\right)^{\alpha} \cdot \varphi
&=& -m\alpha(\alpha-1) \int_{\Omega}\left(u_{\varepsilon}+\varepsilon\right)^{m+\alpha-3}\phi(v_{\varepsilon})\left|\nabla u_{\varepsilon}\right|^{2} \varphi \nn\\
& &   -\alpha(\alpha-1) \int_{\Omega}\left(u_{\varepsilon}+\varepsilon\right)^{m+\alpha-2} \phi'(v_{\varepsilon}) \varphi
        \nabla u_{\varepsilon}\cdot \nabla v_{\varepsilon}  \nn\\
& & -m\alpha\int_{\Omega}\left(u_{\varepsilon}+\varepsilon\right)^{m+\alpha-2} \phi(v_{\varepsilon}) \nabla u_{\varepsilon} \cdot \nabla \varphi \nn\\
& & -\alpha\int_{\Omega} \left(u_{\varepsilon}+\varepsilon\right)^{m+\alpha-1} \phi'(v_{\varepsilon}) \nabla v_{\varepsilon} \cdot \nabla \varphi\nn\\
&=:& I_1 + I_2 + I_3 +I_4
\eea
for all $t>0$ and $\eps\in(0,1)$. Here since Lemma \ref{result3.5} and (\ref{Equ(2.4)}) imply the existence of positive constant $c_2$ such that
\be{uv}
\left\|u_{\varepsilon}(\cdot,t)\right\|_{L^{\infty}(\Omega)} \leq c_2
	\quad \mbox{and} \quad
\left\|\nabla v_{\varepsilon}(\cdot,t)\right\|_{L^{\infty}(\Omega)}
\leq  c_2
\qquad \text {for all $t>0$ and }\eps\in(0,1),
\ee
we have
\be{Equ(3.27)}
I_1
\le \frac{4c_1\kappa_2m \alpha(\alpha-1)}{(m+\alpha-1)^{2}}
\int_{\Omega}\left|\nabla\left(u_{\varepsilon}
+\varepsilon\right)^{\frac{m+\alpha-1}{2}}\right|^{2}
\quad \mbox {for all $t>0$ and $\eps\in(0,1)$}.
\ee
It follows from (\ref{uv}) and the H\"older inequality that
\bea{Equ(3.28)}
I_{2}
&\le& \frac{c_1\alpha(\alpha-1) \kappa_3}{(m+\alpha-1)}
\bigg\{\int_{\Omega}\left|\nabla\left(u_{\varepsilon}+\varepsilon\right)^{m+\alpha-1}\right|^2\bigg\}^{\frac{1}{2}}
\Big\{\io\left|\nabla v_{\varepsilon}\right|^2 \Bigg\}^{\frac{1}{2}}\nn\\
&\le&  \frac{c_1c_2|\Om|^{\frac{1}{2}}\alpha(\alpha-1) \kappa_3}{(m+\alpha-1)}
\bigg\{\int_{\Omega}\left|\nabla\left(u_{\varepsilon}+\varepsilon\right)^{m+\alpha-1}\right|^2\bigg\}^{\frac{1}{2}}
\eea
for all $t>0$ and $\eps\in(0,1)$,
and
\be{Equ(3.29)}
I_{3}
\leq \frac{m\alpha c_1 \kappa_2} {m+\alpha-1}
\bigg\{\int_{\Omega}\left|\nabla\left(u_{\varepsilon}+\varepsilon\right)^{m+\alpha-1}\right|^2 \bigg\}^{\frac{1}{2}}
\quad \mbox {for all $t>0$ and $\eps\in(0,1)$}
\ee
as well as
\bea{Equ(3.30)}
I_{4}
&\le& \kappa_3 \alpha(c_1+1)^{m+\alpha-1}\left\|\nabla v_{\varepsilon}\right\|_{L^{2}(\Omega)}\|\nabla \varphi\|_{L^{2}(\Omega)}\nn\\
&\le&c_1c_2\kappa_3 \alpha(c_1+1)^{m+\alpha-1}
\qquad \mbox {for all $t>0$ and $\eps\in(0,1)$}.
\eea
Collecting (\ref{Equ(3.26)})-(\ref{Equ(3.30)}) and using Young's inequality,
we infer the existence of  $c_3(\alpha)>0$ fulfilling  that
\bas
\left\|\frac{\partial}{\partial t}\left(u_{\varepsilon}+\varepsilon\right)^{\alpha}\right\|_{\left(W_{0}^{k,2}(\Omega)\right)^{*}}
\le  c_3(\alpha)\bigg\{\int_{\Omega}\left|\nabla\left(u_{\varepsilon}+\varepsilon\right)^{\frac{m+\alpha-1}{2}}\right|^{2}
+\int_{\Omega}\left|\nabla\left(u_{\varepsilon}+\varepsilon\right)^{m+\alpha-1}\right|^2
+1\bigg\}
\eas
for all $t>0$ and $\eps\in(0,1)$, and by means of (\ref{Equ(3.21)}) and Lemma \ref{result3.7}, this proves the lemma,
because $\alpha>1$ implies $\frac{m+\alpha-1}{2}>\frac{m}{2}$ and $m+\alpha-1>\frac{m}{2}$.
\qed

Thanks to the boundedness of $u_{\varepsilon}$ and $v_{\varepsilon}$ noted in Lemma \ref{result3.5},
H\"{o}lder estimates for both $v_{\varepsilon}$ and $\nabla v_{\varepsilon}$ can be derived from standard parabolic regularity theory.
\begin{lemma}\label{result3.9}
Let $m > \frac{n}{2}$.  Then there exists $\theta \in(0,1)$ with the property that one can find $C>0$ such that for all $\varepsilon \in(0,1)$,
\begin{equation}\label{Equ(3.33)}
\left\|v_{\varepsilon}\right\|_{C^{\theta,\frac{\theta}{2}}(\bar{\Omega} \times[t, t+1])} \leq C \quad \text { for all } t \geq 0
\end{equation}
and that for all $t_0>0$ it is possible to choose $C(t_0)>0$ fulfilling  that
\begin{equation}\label{Equ(3.34)}
\left\|\nabla v_{\varepsilon}\right\|_{C^{\theta,\frac{\theta}{2}}(\bar{\Omega} \times[t, t+1])} \leq C(t_0) \quad \text { for all } t \geq t_0.
\end{equation}
\end{lemma}
\noindent{\bf{Proof.}} By virtue of the uniform boundedness of $\ueps$ and $\veps$ as shown in Lemma \ref{result3.5},
(\ref{Equ(3.33)}) becomes an immediate consequence of a well-known result on H\"older regularity (\cite{porzio_vespri}),
and (\ref{Equ(3.34)}) can be verified through a reasoning entirely similar to that developed in the proof of \cite[Lemma 3.10]{lankeit_mmm},
which involves maximal Sobolev regularity in parabolic equations and a celebrated embedding result (\cite{amann_2000}),
both applied to a time-localized
version of the second equation in (\ref{Equ(2.1)}). $\hfill{} \Box$

\subsection{Passing to the limit. Proof of Theorem \ref{result1.1}}

With all the preparation in place we can now adopt a standard extraction procedure, leading to the construction of a limit function $(u,v)$, and moreover,
such a limit function is proved to be a weak solution to (\ref{Equ(1.1)}) as documented in Theorem \ref{result1.1}.

\begin{lemma}\label{result3.10}
Let $m > \frac{n}{2}$. Then there exist $\left(\varepsilon_{j}\right)_{j \in \mathbb{N}} \subset(0,1)$ and  functions $u$ and $v$ such that $\varepsilon_{j} \searrow 0$ as $j \rightarrow \infty$, and that
\begin{eqnarray}
& & u_{\varepsilon}\rightarrow u
\qquad \text {a.e.~in  $\Omega\times (0, \infty)$ and in $L_{\text {loc }}^{p}(\bom\times [0,\infty))$ for all $p\in [1,\infty)$},
\label{Equ(3.101)}\\
& &  u_{\varepsilon} \stackrel{\star}{\rightharpoonup} u
\qquad \text {in } L^{\infty}(\Omega \times(0, \infty)),
\label{Equ(3.102)}\\
& & \na (\ueps+\eps)^m \wto \na \ueps^m
\qquad \mbox{in } L^2_{loc}(\bom\times [0,\infty)),
\label{3103}\\
%& & u_{\varepsilon} \rightarrow u \text { in } \bigcap_{p \geq 1} L_{\text {loc }}^{p}(\bar{\Omega} \times[0, \infty))
%\label{Equ(3.103)}\\
& & v_{\varepsilon} \rightarrow v
\qquad \text {in } C_{\text {loc }}^{0}(\bar{\Omega} \times[0, \infty)),
\label{Equ(3.104)}\\
& & \nabla v_{\varepsilon} \rightarrow \nabla v
\qquad \text {in } C_{\text {loc }}^{0}(\bar{\Omega} \times((0, \infty)),
\label{Equ(3.105)}
%& & \phi\left(v_{\varepsilon}\right) \rightarrow \phi(v)
%\qquad \text {in } C_{\text {loc }}^{0}(\bar{\Omega} \times([0, \infty)),
%\label{Equ(3.1061)}
%& & v_{\varepsilon} \stackrel{\star}{\rightharpoonup} v
%\qquad \text {in } L^{\infty}\left(\Omega\times(0, \infty) \right)
%\label{Equ(3.106)}
\end{eqnarray}
as $\varepsilon=\varepsilon_{j} \searrow 0.$ Moreover, $(u,v)$
forms a global weak solution of (\ref{Equ(1.1)})  in the sense of Definition \ref{result0.1}, and we have
\begin{equation}\label{3107}
\int_{\Omega} u(\cdot, t)=\int_{\Omega} u_{0} \quad \text { for a.e. } t>0.
\end{equation}
\end{lemma}
\noindent{\bf{Proof.}} Since $n\ge2$ and $m>\max\{\frac{m}{2}, \, \frac{n}{2}\}$, then fixing $k\in\mathbb{N}$
such that $k>\frac{n}{2}$, we know from Lemmata \ref{result3.5}, \ref{result3.7} and \ref{result3.8}
that for all $T>0$,
 \be{31001}
  ((\ueps+\eps)^m)_{\eps\in(0,1)}
  \mbox{ is bounded in $L^2((0,T); W^{1,2}(\Om))$}
 \ee
 and
 \be{31002}
  (\pa_t(\ueps+\eps)^m)_{\eps\in(0,1)}
  \mbox{ is bounded in $L^1((0,T); (W_0^{k,2}(\Om))^{\star})$},
 \ee
 so that, by (\ref{31001}) and (\ref{31002}), we infer from an  Aubin-Lions lemma (\cite{temam}) and the compact embedding
 $W^{1,2}(\Om) \hra L^2(\Om)$  that there exist $(\eps_j)_{j\in\mathbb{N}}\subset(0,1)$ and a nonnegative
 function $u$ defined on $\Om\times (0,\infty)$ such that $\eps=\eps_j\searrow 0$ as $j\to\infty$, and that
 \be{31003}
  (\ueps+\eps)^m \rightarrow u^m
  \qquad\mbox{in } L^2_{loc}(\bom\times[0,\infty)),
 \ee
as $\eps=\eps_j\searrow 0$, hence (\ref{3103}) holds and
 \be{31004}
  u_{\varepsilon}\rightarrow u
  \qquad\mbox{in $L^{2m}_{loc}(\bom\times[0,\infty))$ and a.e. in $\Om\times(0,\infty)$},
 \ee
as $\eps=\eps_j\searrow 0$. By Lemma \ref{result3.5} and Egorov's theorem, we infer from (\ref{31004}) that
(\ref{Equ(3.102)}) and (\ref{Equ(3.101)}) hold. Also it follows from (\ref{31004}) that (\ref{3107}) is true
by the Fubini-Tonelli theorem and (\ref{Equ(2.2)}).
Next upon two applications of the Arzel\`a-Ascoli theorem, we can readily verify (\ref{Equ(3.104)}) and (\ref{Equ(3.105)}) with
$v\in L^{\infty}(\Om\times(0,\infty))\cap C_{\text {loc }}^{0}(\bar{\Omega} \times[0, \infty))$ and $\nabla v\in
C_{\text {loc }}^{0}(\bar{\Omega} \times((0, \infty))$, and in particular, by (\ref{Equ(3.104)}), (\ref{phi}) and (\ref{Equ(2.3)}),
we deduce that
 \be{31005}
  \phi(\veps) \rightarrow  \phi(v)
  \qquad\mbox{in $C_{\text {loc }}^{0}(\bar{\Omega} \times[0, \infty))$}
 \ee
 and
 \be{31006}
 \phi'(\veps) \rightarrow  \phi'(v)
   \qquad\mbox{in $C_{\text {loc }}^{0}(\bar{\Omega} \times[0, \infty))$}
 \ee
 as $\eps=\eps_j\searrow 0$. On the basis of (\ref{31006}), (\ref{Equ(3.105)}), (\ref{phi}), and (\ref{Equ(2.3)}),
 we have
  \be{phi'na}
   \phi'(\veps)\na \veps \rightarrow \phi'(v)\na v
   \qquad\mbox{in $L^2_{loc}(\bom\times[0,\infty))$},
  \ee
  as $\eps=\eps_j\searrow 0$. \abs
  Now for fixed $\vp\in C^{\infty}_0(\bom\times[0,\infty))$,
  from (\ref{Equ(2.1)}) it follows that
   \bea{ue}
  	- \int_0^\infty \io \ueps\vp_t - \io u_0\vp(\cdot,0)
 &=& - \int_0^\infty \io\phi(\veps)\na (\ueps+\eps)^m\cdot\na \vp\nn\\
  & & - \int_0^\infty \io\phi'(\veps) (\ueps+\eps)^m\na \veps\cdot\na \vp
    \eea
for all $\eps\in(0,1)$. In view of (\ref{Equ(3.101)})-(\ref{Equ(3.105)}) and (\ref{31005})-(\ref{phi'na}), taking $\eps=\eps_j\searrow 0$
we have (\ref{Equ(1.6)}). The verification of (\ref{Equ(1.7)}) follows in a similar manner, and moreover, the regularity requirements
(\ref{Equ(1.5)}) and (\ref{w1}) recorded in Definition \ref{result0.1} become straightforward consequences of (\ref{Equ(3.101)})-(\ref{Equ(3.105)}).
$\hfill{} \Box$

We are now in a position to verify our main result.\\
{\bf Proof of Theorem \ref{result1.1}.}  Theorem \ref{result1.1} is a direct consequence of Lemmata \ref{result3.5} and \ref{result3.10}.\qed

{\bf Acknowledgments}
The first author was funded by the China Scholarship Council (No. 202006630070).
The second author was supported by the China Scholarship Council (No. 202108500085) and Natural Science Foundation of Chongqing (No. cstc2021jcyj-msxmX0412).

\end{document}